\documentclass{article}
\usepackage{graphicx}
\usepackage{epsfig}
\usepackage{amssymb, amsmath}
\usepackage{amsthm}
\usepackage{setspace}
\usepackage[top = 1in, bottom = 1in, left = 1.2in, right =
1.2in]{geometry}
\newtheorem{theorem}{Theorem}[section]

\newtheorem{lemma}[theorem]{Lemma}

\newtheorem{remark}[theorem]{Remark}

\begin{document}

\title{The 2D Euler-Poisson System with Spherical Symmetry}%
\author{Juhi Jang\\
University of California Riverside\\
juhijang@math.ucr.edu}%

\singlespacing \maketitle \numberwithin{equation}{section}

\begin{abstract}
This article concerns the global-in-time existence of smooth solutions with small amplitude to two space dimensional Euler-Poisson system. The main difficulty lies in the slow time decay $(1+t)^{-1}$ of the linear system. Inspired by Ozawa, Tsutaya, and Tsutsumi's work \cite{ott,ott2}, we show that such smooth solutions exist for radially symmetric flows. 
\end{abstract}


\maketitle

\section{Introduction and Main result}

More than 99$\%$ of the state of matter in the universe is plasma, a
collection of fast-moving dilute electrons and ions. Plasma also
appears in nuclear fusion, semiconductors, optics, computers, and so
on. The hydrodynamic model for plasma treats the propagation of
electrons as the flow of compressible charged fluid. One of the
simplest hydrodynamic models is the Euler-Poisson system in an
electric field with the repulsive Coulomb interaction:
\begin{equation}
\begin{split}
\partial_t n+\nabla\cdot(n u)&=0,\\
\partial_t u+u\cdot\nabla u+\frac{1}{m_en}\nabla p(n)&=
\frac{e}{m_e}\nabla\phi,\\
\Delta\phi=4\pi e(n-n_0),\;|\phi|\rightarrow 0, &\;\text{as }|x|
\rightarrow \infty,\label{epp}
\end{split}
\end{equation}
where the electrons of charge $e$ and mass $m_e$ are described by a
density $n(t,x)$, an average velocity $u(t,x)$, and 
 $\nabla\phi(t,x)$ is the electric field for $t\geq 0$ and $x\in\mathbb{R}^d$, and $en_0$
represents the constant equilibrium charge density. The pressure is
assumed to obey the $\gamma$-law; that is, $p(n)=An^{\gamma}$ for
$\gamma\geq 1$ and an entropy constant $A$. 

The purpose of this article is to investigate whether the system
(\ref{epp}) when the spatial dimension is given by $d=2$ admits global smooth solutions in time for small
perturbed initial data. 

In the absence of the electric field
($\nabla\phi\equiv 0$), (\ref{epp}) leads to the compressible Euler
equations, and it is well known that finite time singularity, as
shock waves, develops generically even for small data \cite{si}. On
the other hand, the Euler-Poisson system (\ref{epp}) displays rather
contrasting phenomena due to the dispersive effect of Poisson
forcing: the system can be written as  quasilinear Klein-Gordon equations with quadratic nonlinearity under the invariant irrotationality  condition $\nabla \times u=0$. In order to see this, following \cite{g1}, we first define the following new variables:
\begin{equation}
m(t,x)\equiv
\frac{2}{\gamma-1}\{(\frac{n(\frac{t}{c_0},x)}{n_0})^{\frac{\gamma-1}{2}}-1\},\;
\;v(t,x)\equiv \frac{1}{c_0}u(\frac{t}{c_0},x),\;\;\psi(t,x)\equiv
\phi(\frac{t}{c_0},x)
\end{equation}
where $c_0=\sqrt{\gamma}n_0^{\frac{\gamma-1}{2}}$ is the sound
speed. In terms of new variables, the Euler-Poisson system
(\ref{epp}) takes the form of
\begin{equation}
\begin{split}
&\partial_tm+\nabla\cdot v+v\cdot\nabla
m+\frac{\gamma-1}{2}m\nabla\cdot v=0\\
&\partial_tv+\nabla m+v\cdot \nabla v +\frac{\gamma-1}{2}m\nabla m
=\frac{1}{c_0^2}\nabla \psi\\
&\partial_t\nabla\psi=-n_0v-n_0\nabla\Delta^{-1}\nabla\cdot\{[(\frac{\gamma-1}{2}
m+1)^{\frac{2}{\gamma-1}}-1]v\}\label{ep}
\end{split}
\end{equation}
with the constraint
\[
\Delta\psi=n_0[(\frac{\gamma-1}{2}
m+1)^{\frac{2}{\gamma-1}}-1]\equiv n_0[m-h(m)].
\]
To derive the last equation for $\psi$ in (\ref{ep}), we have used
the irrotationality condition $\nabla \times u=0$. Now by taking one
more derivative of (\ref{ep}) we obtain
\begin{equation}
\begin{split}
 (\partial_t^2-\Delta+m_0)m=&\nabla\cdot[v\cdot\nabla v+\frac{\gamma-1}{2} m\nabla m]\\
&-\partial_t[v\cdot\nabla m+\frac{\gamma-1}{2}m\nabla\cdot v]+m_0h(m)\\
(\partial_t^2-\Delta+m_0)v=&\nabla[v\cdot\nabla m+\frac{\gamma-1}{2}m\nabla\cdot v] -\partial_t[v\cdot\nabla v+\frac{\gamma-1}{2} m\nabla m]\\
&-m_0\nabla\Delta^{-1}\nabla\cdot\{[m-h(m)]v\}\label{kg}
\end{split}
\end{equation}
where $m_0=\frac{n_0}{c_0^2}$. Note that the right-hand-side in
(\ref{kg}) is formally second order in $m$ and $v$, and therefore, the new equations 
(\ref{kg}) are the \textit{quasilinear Klein-Gordon equations with quadratic
nonlinearity} and moreover, the nonlinearity contains the \textit{non-local} term
$m_0\nabla\Delta^{-1}\nabla\cdot\{[m-h(m)]v\}$.  Note that $\nabla\Delta^{-1}\nabla$ is an example of singular integrals, so called the Riesz potential. 

In \cite{g1}, Guo constructed global-in-time smooth
irrotational flows with small perturbed initial data in
$\mathbb{R}^3$, based on the above Klein-Gordon effect that gives a time
decay rate of $(1+t)^{-\frac{3}{2}}$, by making use of Shatah's
normal form method \cite{s} and also by using the theory of singular integrals to deal with non-local terms. However, the analysis there is not
directly applicable to 2D flows due to the slower time decay rate
$(1+t)^{-1}$, which is not integrable in time, of the 2D linear Klein-Gordon equation. In \cite{ott,ott2},
Ozawa, Tsutaya, and Tsutsumi overcame this difficulty  by combining
the normal form method of Shatah  \cite{s} and the decay
estimate of the linear Klein-Gordon equation due to Georgiev \cite{g} based on the vector field method of
Klainerman \cite{k}, and they proved the global existence of 2D
quasilinear Klein-Gordon equations with quadratic (but local) nonlinearity with small amplitude in high-order Sobolev spaces. In particular, their results imply the existence of $C^\infty$ solutions. Inspired by this result, we establish the following global-in-time existence of the Euler-Poisson system:

\begin{theorem}There exists unique global
smooth solutions to the 2D Euler-Poisson system (\ref{epp}) for
spherically symmetric flows with small perturbed initial data.\label{thm}
\end{theorem}

What follows next are the brief sketch of the proof of Theorem \ref{thm} and other recent progress on the Euler-Poisson system. 

\section{Proof of Theorem \ref{thm}}

Note that from the first and third equation in \eqref{epp} the dynamics of the electric field $\nabla\phi$ is determined by 
\begin{equation}\label{phi}
\partial_t \nabla\phi = -\nabla \triangle ^{-1} \nabla \cdot[nu]
\end{equation}
and also see the third equation in \eqref{ep}. In general, this non-local interaction makes the analysis intriguing.
Nevertheless, we can show that the non-local effect is not present for radially symmetric irrotational flows. To see that, we start with the following two elementary lemmas. 

\begin{lemma}
 Let $\eta$ be a vector field with suitable regularity in
$\mathbb{R}^{3}$ (or $\mathbb{R}^{2}$). If $\nabla \times \eta = 0$,
then $\nabla \triangle ^{-1} \nabla \eta = \eta$. 
\end{lemma}

\begin{proof}
Since $\mathbb{R}^{3}$(or $\mathbb{R}^{2}$) is simply
connected, the condition $\nabla \times \eta = 0$ implies that $\eta$ is a
gradient field. Hence, there exists a scalar function  in $\mathbb{R}^{3}$(or $\mathbb{R}^{2}$) such that $\eta = \nabla \varphi$. By taking the
divergence of both sides, we see that $\nabla\cdot \eta = \triangle \varphi$. Hence,
$\triangle ^{-1}  \nabla \cdot\eta = \varphi$. Next by 
taking the gradient, we get the desired result.
\end{proof}

The next lemma shows that radially symmetric vector fields in $\mathbb{R}^2$ automatically satisfy the curl-free condition,  namely, the irrotationality condition. For the purpose of the application to our problem, we state the result in terms of $nu$ in \eqref{phi}. 

\begin{lemma} Let $n$ be a radially symmetric function and let $u$ be a radially symmetric vector field in $\mathbb{R}^2$ so that $n(x)=n(|x|)$ and  $u(x)=u(|x|)\frac{x}{|x|}$.    Then the following holds: 
 $$\nabla \times [nu] = 0.$$ 
\end{lemma}

\begin{proof} Let $r=|x|$. Using the polar coordinates $(r,\theta)$, write
$nu=(\cos\theta,\sin\theta)\varphi(r)$ for a scalar function
$\varphi$. Note $\partial_x = \cos\theta\partial_r -
{\sin\theta\over r}\partial_\theta$ and $\partial_y =
\sin\theta\partial_r + {\cos\theta\over r}\partial_\theta.$ Thus we obtain 
\begin{equation*}
\begin{split}\nabla\times [nu] &=
(\cos\theta\partial_r - {\sin\theta\over r}\partial_\theta)
\sin\theta\varphi(r) -(\sin\theta\partial_r + {\cos\theta\over
r}\partial_\theta)
\cos\theta\varphi(r)\\
&=[\cos\theta \sin\theta \varphi'(r)- {\sin\theta\over
r}\cos\theta\varphi(r)]-[\sin\theta\cos\theta\varphi'(r)+{\cos\theta\over
r}(-\sin\theta)\varphi(r)]\\
&=0.
\end{split}
\end{equation*}\end{proof}

Hence, by the above two lemmas, we see that for radially symmetric flows, the dynamics of $\nabla\phi$ in \eqref{phi} is given by $$\partial_t \nabla\phi = - nu.$$
This fact also confirms that the non-local term $m_0\nabla\Delta^{-1}\nabla\cdot\{[m-h(m)]v\}$ in \eqref{kg} becomes local. Hence, we are now ready to apply the global existence result of Ozawa, Tsutaya, and Tsutsumi \cite{ott,ott2} to \eqref{kg}, the quasilinear Klein-Gordon equations with quadratic nonlinearity  for two space dimension and therefore, accordingly we obtain the global-in-time existence of spherically symmetric  
flows to the Euler-Poisson System in $\mathbb{R}^{2+1}$. We refer to \cite{ott,ott2} for the  function spaces in which the solutions $m,v$ to \eqref{kg} reside. This finishes the proof of Theorem \ref{thm}.

\section{Discussion}

Some remarks follow in the below. 

\begin{remark}
 Lemma 2.2 is not valid for general irrotational flows and thus Theorem \ref{thm} may not be directly extended to the general irrotational flows. For general 2D flow, the difficulty comes from not only the slow dispersion of the linear part  but also the non-local quadratic nonlinearity  in \eqref{kg}. Up to the author's knowledge, the global existence involving the non-local feature of the Cauchy problem to \eqref{kg} even for irrotational flows is still open. To tackle the 2D Cauchy problem, an advanced sharp technique with aid of harmonic analysis seems to be needed. 
\end{remark}

\begin{remark}
 On the other hand, recently, in \cite{jlz}  by Li, Zhang, and the author, some global solutions to 2D Euler-Poisson system  were constructed by solving the final data problem, in other words by  constructing the wave operators to \eqref{kg}. 
\end{remark}

\begin{remark}
 For 1D flows, in \cite{tw},
Tadmor and Wei showed that (\ref{epp}) admits a global smooth
solution for certain initial configurations by clarifying an
intrinsic critical threshold. However, in multi-dimensional flows,
the critical threshold phenomenon has been captured only for
pressureless restricted Euler-Poisson models \cite{elt}.
\end{remark}

\begin{remark}  Recently,   Guo and Pausader \cite{gp} have established the global existence  to the 3D Euler-Poisson  system describing the Ion dynamics. The question is open for the Euler-Poisson system of the full electron and ion dynamics of two-fluid for both two and three dimensional cases. 
\end{remark}

\

\textbf{Acknowledgments:} The author would like to thank \textsc{Yan Guo} to encourage her to write this note.

\end{document}